\documentclass[12pt]{article}
\usepackage{amsmath,amsthm}
\usepackage{amssymb,latexsym}
\usepackage{enumerate}
\usepackage[english]{babel}
\usepackage{graphicx}

\def\Q{{\mathbb Q}}
\def\Z{{\mathbb Z}}

\def\Re{{\rm Re}}
\def\Im{{\rm Im}}

\newtheorem{lemma}{Lemma}
\newtheorem{theorem}[lemma]{Theorem}

\title{
Calculating "small" solutions of  \\ relative Thue equations
}
\author{
Istv\'{a}n Ga\'{a}l\thanks{
        This work was partially supported by the European Union 
        and the European Social Fund through project Supercomputer, 
        the national virtual lab (grant no.: TAMOP-4.2.2.C-11/1/KONV-2012-0010)
        and also supported in part by K67580 and K75566 from the
        Hungarian National Foundation for Scientific Research,
         }
\\
University of Debrecen, Mathematical Institute \\
            H--4010 Debrecen Pf.12., Hungary \\
            e--mail: igaal@science.unideb.hu\\ \\
}

\begin{document}

\maketitle
\thispagestyle{empty}

\renewcommand{\thefootnote}{}

\footnote{2010 \emph{Mathematics Subject Classification}: Primary 11D59, 11Y50; Secondary 11R21}

\footnote{\emph{Key words and phrases}: relative Thue equations; 
algoriths for the resolution; LLL reduction; supercomputers}

\renewcommand{\thefootnote}{\arabic{footnote}}
\setcounter{footnote}{0}

\begin{abstract}
Diophantine equations can often be reduced to various types of 
classical Thue equations \cite{thue}, \cite{baker}. 
These equations usually have only
very small solutions, on the other hand to compute all solutions (i.e. to 
prove the non-existence of large solutions) is a time consuming procedure.
Therefore it is very practical to have a fast algorithm to calculate the "small"
solutions, especially if "small" means less than e.g. $10^{100}$. Such an
algorithm was constructed by A.Peth\H o \cite{pet} in 1987 based
on continued fractions.

In the present paper we construct a similar type of fast algorithm 
to calculate "small" solutions of {\it relative} Thue equations. 
Our method is based on the LLL reduction algorithm.
We illustrate the method with explicit examples.
The algorithm has several applications.
\end{abstract}

\section{Introduction}

Let $F\in \Z[x,y]$ be a binary form of degree $\geq 3$, irreducible over $\Q$.
There is an extensive literature (cf. \cite{thue}, \cite{baker},
\cite{book}) on {\it Thue equations} of the type
\begin{equation}
F(x,y)=m \;\; {\rm in}\; x,y\in\Z.
\label{ttt}
\end{equation}
Various types of diophantine equations, among others index form equations
(cf. \cite{thueappl}, \cite{gppsim}) can be reduced to Thue equations.
Using the effective method of A.Baker \cite{baker} and reduction methods
(for a survey see \cite{book}) there is an algorithm for solving Thue
equations completely. However, this procedure is rather time consuming
while our experience shows that such equations only have small solutions.
Hence the efforts are invested in the proof of the nonexistence
of large solutions, not in the calculation of the solutions.

Therefore in many applications and practical methods 
the fast algorithm of
A.Peth\H o \cite{pet}, giving only the solutions with $|y|<C$ was useful,
especially because it remains fast also for e.g. $C=10^{100}$.
(In \cite{purecubic} we calculated "small" solutions of thousands of
index form equations  in pure cubic fields, in \cite{gr} we calculated "small"
solutions of binomial Thue equations $x^4-my^4=\pm 1$ for $0\leq m\leq 10^7$.)

Note that if $F$ has leading coefficient 1 (which can be assumed without
restricting the generality), $\alpha$ is a root of $F(x,1)=0$ 
and $K=\Q(\alpha)$, then equation (\ref{ttt})
can be written in the form 
\begin{equation}
N_{K/Q}(x-\alpha y)=m \;\; {\rm in}\; x,y\in\Z.
\label{cc}
\end{equation}

Let $M$ be an algebraic number field with ring of integers $\Z_M$.
Let $\alpha$ be an algebraic integer over $M$ and set $K=M(\alpha)$. 
Let $0\neq \mu\in \Z_M$ and consider the {\it relative Thue equation}
\begin{equation}
N_{K/M}(X-\alpha Y)=\mu \;\;\; {\rm in} \;\;\; X,Y\in\Z_M.
\label{rel}
\end{equation}
This equation is a direct analogue of (\ref{cc}).

Diophantine problems often lead also to relative Thue
equations. The index form equations of sextic fields with a quadratic subfield
\cite{c6}, \cite{csextic}, of nonic fields with a cubic subfield \cite{nine}, 
of quartic relative extensions \cite{g42}, all lead to relative Thue equations.

There is an algorithm for the complete resolution of relative Thue equations by
I.Ga\'al and M.Pohst \cite{relthue}. This involves Baker's method, a reduction
procedure using LLL and an enumeration method (see also \cite{book}).
The execution time takes some hours.
As initial data the method needs the fundamental units of the number field
involved and calculation of all non-associated elements of given relative norm.
Note that often the main difficulty is to calculate these basic data
in a higher degree field (of degree $\geq 10$), requiring again considerable
CPU time. For fields of degree $\geq 15$ these procedures often fail.

It is therefore useful to have a faster algorithm 
calculating the "small" solutions of relative Thue equations
that does not require initial data (e.g. fundamental units, elements 
of given relative norm) of the number field and 
is applicable also for higher degree relative Thue equations. In this paper
our purpose is to construct such an algorithm.

As we shall see our new algorithm works efficiently even for higher degree relative
Thue equations that cannot be attacted by the methods of \cite{relthue}.
It is very efficient over quadratic number fields $M$ 
but usable also over cubic and quartic fields.

\section{Elementary estimates for relative Thue equations}
\label{sectel}

Let $m=[M:\Q]$ and $(1,\omega_2,\ldots,\omega_m)$ be an integral basis of the ring 
of integers of $M$.
Denote by $\gamma^{(j)},\; (j=1,\ldots, m)$ the conjugates of any $\gamma\in M$.

Denote by $n$ the degree of $\alpha$ over $M$ and by $f(x)$ its relative defining polynomial over $M$. Denote by $\alpha^{(jk)}, \; (k=1,\ldots, n)$ the relative conjugates 
of $\alpha$ over $M^{(j)}$, that is the roots of the $j$-th conjugate of $f(x)$.
We also denote by $\gamma^{(jk)}$ the conjugates of any $\gamma\in K$ corresponding to
$\alpha^{(jk)}$.

We assume that $n>m$ if $K$ is not totally real and $n>2m$ if $K$ is totally real.
This condition ensures that our reduction procedure in Section \ref{sectred}
is efficient.

Our purpose is to determine all solutions of (\ref{rel})
with $\overline{|Y|}<C$, where $C$ is a large
given constant, say $10^{100}$ or $10^{500}$.
(For any algebraic number $\gamma$ we denote by $\overline{|\gamma|}$ the {\it size} 
of $\gamma$, that is the maximum absolute value of its conjugates.)

We represent $X$ and $Y$ in the form
\[
X=x_1+\omega_2x_2+\ldots +\omega_m x_m,\;\;\; Y=y_1+\omega_2y_2+\ldots +\omega_m y_m,
\]
with $x_i,y_i\in\Z \; (1\leq i\leq m)$. We set $A=\max (\max |x_i|,\max |y_i|)$. 

Our algorithm is based on the well known fact that for any fixed solution $X,Y\in \Z_M$
of equation (\ref{rel}) there is a conjugate of $\beta=X-\alpha Y$ which is very small, assumed that
$Y$ is not very small.

To formulate the main result we need the following notation.

Let 
\[
c_{hji}=\frac{1}{2}\left|(\alpha^{(hj)}-\alpha^{(hi)})\right|,\;\; {\rm for} \;
1\leq h\leq m,\;\;1\leq i,j\leq n, \;\;i\neq j,
\]
\[
c_{hi}=\frac{|\mu^{(h)}|}{ \prod_{1\leq j\leq n, j\neq i} c_{hji} }, \;\; {\rm for} \;
1\leq h\leq m,\;\;1\leq i\leq n, 
\]
\[
c_1=\max_{h,j,i}  \frac{\sqrt[n]{|\mu^{(h)}|}}{c_{hji}}, 
\]
where the maximum is taken for $1\leq h\leq m,\;\;1\leq i,j\leq n, \;\;i\neq j$.

Let $S$ be the $m\times m$ matrix with entries $1,\omega_2^{(j)},\ldots,\omega_m^{(j)}$
in the $j$-th row. Denote by $c_2$ the row norm of $S^{-1}$ that is the maximum sum
of the absolute values of the elements in its rows.

Let
\[
c_3=\frac{\sqrt[n]{\overline{|\mu|}}}{\overline{|\alpha|}},
\]
\[
c_4=\max (c_1,c_3),
\]
\[
c_5=2 c_2 {\overline{|\alpha|}} 
\]
and
\[
d_{hi}=c_{hi} c_5^{n-1}\;\; {\rm for} \;
1\leq h\leq m,\;\;1\leq i\leq n.
\]
Finally, let $C$ be a given constant. 

\begin{theorem}
If $(X,Y)\in\Z_M^2$ is a solution of equation (\ref{rel}) with $\overline{|Y|}>c_3$, then 
\begin{equation}
A\leq c_5 \cdot \overline{|Y|}.
\label{aa}
\end{equation}
 \label{th1}
\end{theorem}

\noindent
{\bf Proof of Theorem \ref{th1}}.\\
Set $\beta=X-\alpha Y$.
For any $k$ let $\ell$ be the index with $|\beta^{(k\ell)}|=\min_{1\leq j\leq n}|\beta^{(kj)}|$.
Equation (\ref{rel}) implies
\[
\beta^{(k1)}\ldots \beta^{(kn)}=\mu^{(k)}.
\]
Therefore $|\beta^{(k\ell)}|\leq \sqrt[n]{|\mu^{(k)}|}$ whence we have
\[
|X^{(k)}|\leq |\beta^{(k\ell)}| + |\alpha^{(k\ell)}|\cdot |Y^{(k)}|
\leq \sqrt[n]{\overline{|\mu|}} + \overline{|\alpha|} \cdot \overline{|Y|}
\]
whence using $\overline{|Y|}>c_3$ we obtain
\begin{equation}
\overline{|X|}\leq 2 \; \overline{|\alpha|} \cdot \overline{|Y|}.
\label{xx}
\end{equation}
By
\[
\left(\begin{array}{c}y_1\\  \vdots   \\  y_n \end{array}  \right)
=
S^{-1}
\left(\begin{array}{c}Y^{(1)}\\  \vdots   \\  Y^{(m)} \end{array}  \right)
\]
we obtain $\max|y_j|\leq c_2 \overline{|Y|}$. Similarly we have $\max|x_j|\leq c_2 \overline{|X|}$,
whence by (\ref{xx}) we get the assertion (\ref{aa}).
\hfill $\Box$\\

\vspace{0.5cm}

\begin{theorem}
If $(X,Y)\in\Z_M^2$ is a solution of equation (\ref{rel}) with
$\overline{|Y|}> c_4$, then there exist $h,i$ ($1\leq h\leq m,\;\;1\leq i\leq n$),
such that
\begin{equation}
|\beta^{(hi)}|\leq d_{hi} A^{1-n}.
\label{bb}
\end{equation}
 \label{th2}
\end{theorem}

\noindent
{\bf Proof of Theorem \ref{th2}}.\\
Let $Y^{(h)}$ be the conjugate of $Y$ with $\overline{|Y|}=|Y^{(h)}|$ and let $i$
be determined by
\[
|\beta^{(hi)}|=\min_{1\leq j\leq n}|\beta^{(hj)}|.
\]
Obviously
\begin{equation}
|\beta^{(hi)}|\leq \sqrt[n]{|\mu^{(h)}|}
\label{kicsi}
\end{equation}
and for any $j\neq i \;\; (1\leq j\leq n)$ using $\overline{|Y|}>c_1$ we have
\begin{equation}
|\beta^{(hj)}|\geq   |\beta^{(hj)}-\beta^{(hi)}|-|\beta^{(hi)}|
              \geq   |(\alpha^{(hj)}-\alpha^{(hi)})Y^{(h)}|-\sqrt[n]{|\mu^{(h)}|}
              \geq c_{hji} \overline{|Y|}
\label{tobbi}
\end{equation}
By equation (\ref{rel}) we have now
\[
\beta^{(h1)}\ldots \beta^{(hn)}=\mu^{(h)}
\]
therefore
\begin{equation}
|\beta^{(hi)}|\leq c_{hi} \cdot \overline{|Y|}^{\; 1-n}
\label{kisebb}
\end{equation}
By $\overline{|Y|}>c_3$ Theorem \ref{th1} applies, therefore (\ref{aa}) is satisfied, whence
${\overline{|Y|}\;}^{-1}\leq c_5\;A^{-1}$, that is 
\begin{equation}
{\overline{|Y|}\;}^{1-n}\leq c_5^{n-1} A^{1-n}.
\label{aa2}
\end{equation}
Therefore by inequality (\ref{kisebb}) we obtain the assertion (\ref{bb}).
\hfill $\Box$\\

\vspace{0.5cm}

\section{Reducing the bound for $A$}
\label{sectred}

In this section we develop a reduction procedure for $A$.
More exactly we apply the extension of M.Pohst \cite{dmv} of the standard LLL algorithm
of A.K.Lenstra, H.W.Lenstra Jr. and L.Lov\'asz \cite{lll}.

Let $H$ be a large constant to be given later, let $i,h$
be given indices with $1\leq h\leq m,\;\;1\leq i\leq n$ such that
inequality (\ref{bb}) is satisfied, that is 
\[
|x_1+\omega_2^{(h)}x_2+\ldots + \omega_m^{(h)}x_m 
-\alpha^{(hi)}y_1-\alpha^{(hi)}\omega_2^{(h)}y_2-\ldots - \alpha^{(hi)}\omega_m^{(h)}y_m | \leq 
\]
\begin{equation}
\leq d_{hi} A^{1-n}.
\label{redineq}
\end{equation}

Consider now the lattice $\cal{L}$ generated by the columns of the matrix
{\small
\[
\cal{L}=
\left(
\begin{array}{ccccccccc}
1&0&\ldots&0&\vline&0&0&\ldots &0\\
0&1&\ldots&0&\vline&0&0&\ldots &0\\
\vdots&\vdots&\ddots&\vdots&\vline&\vdots&\vdots&\ddots&\vdots\\
0&0&\ldots&1&\vline&0&0&\ldots &0\\
\hline
0&0&\ldots&0&\vline&1&0&\ldots &0\\
0&0&\ldots&0&\vline&0&1&\ldots &0\\
\vdots&\vdots&\ddots&\vdots&\vline&\vdots&\vdots&\ddots&\vdots\\
0&0&\ldots&0&\vline&0&0&\ldots &1\\
\hline
H&H\Re(\omega_2^{(h)})&\ldots & H\Re(\omega_m^{(h)})&\vline&
H\Re(\alpha^{(hi)})&H\Re(\alpha^{(hi)}\omega_2^{(h)})&\ldots & H\Re(\alpha^{(hi)}\omega_m^{(h)})\\
H&H\Im(\omega_2^{(h)})&\ldots & H\Im(\omega_m^{(h)})&\vline&
H\Im(\alpha^{(hi)})&H\Im(\alpha^{(hi)}\omega_2^{(h)})&\ldots & H\Im(\alpha^{(hi)}\omega_m^{(h)})\\
\end{array}
\right)
\]
}
In the totally real case we may omit the last row.
Denote by $b_1$ the first vector of an LLL-reduced basis of the lattice ${\cal L}$.

\begin{theorem}
\label{redlemma}
Assume that $x_1,\ldots,x_m,y_1,\ldots,y_m$ are integers 
with $A=\max(\max|x_i|,\max|y_i|)$, such that (\ref{redineq}) is satisfied.
If $A\leq A_0$ for some constant $A_0$ and for the first vector $b_1$
of the LLL reduced basis of $\cal L$ we have
\[
|b_1|\geq \sqrt{(2m+1)2^{2m-1}} \cdot A_0,
\]
then
\begin{equation}
A\leq \left( \frac{d_{hi} H}{A_0}\right)^{\frac{1}{n-1}}
\label{rrr}
\end{equation}
\label{th3}
\end{theorem}

\noindent
{\bf Proof of Theorem \ref{th3}}.\\ 
The proof is almost the same as in \cite{book}.
Denote by $l_0$ the shortest vector in the lattice $\cal L$.
Assume that the vector $l$ is a linear combination of the
lattice vectors with integer coefficients $x_1,\ldots,x_m,y_1,\ldots,y_m$,
respectively. Observe that the last two components of $l$ are 
the real and imaginary parts of $\beta^{(hi)}$.

Using the inequalities of \cite{dmv} we have $|b_1|^2\leq 2^{2m-1} |l_0|^2$. 
Obviously $|l_0|\leq |l|$.
The first $2m$ components of $l$ are in absolute value $\leq A_0$,
for the last 2 components (\ref{redineq}) is satisfied. 
Hence we obtain
\begin{eqnarray*}
&&(2m+1)A_0^2 = 2^{1-2m}\left( (2m+1)\cdot 2^{2m-1} A_0^2 \right)  \\
&&\leq 2^{1-2m}|b_1|^2\leq |l_0|^2\leq |l|^2
\leq 2m \cdot A_0^2+H^2d_{hi}^2 A^{2-2n},
\end{eqnarray*}
whence 
\[
A_0\leq d_{hi} H A^{1-n} \;\; ,
\]
which implies the assertion.
\hfill $\Box$\\

\vspace{0.5cm}

\noindent
{\bf Remark 1}\\
We made several tests to figure out how one can suitably choose $H$, for which $m$ and $n$
the procedure is applicable and what is the magnitude of the reduced bound.
We summarize our experiences in the following table:
\[
\begin{array}{|c|c|c|}
\hline
                                 & {\rm complex}\;{\rm case} & {\rm totally}\;{\rm real}\;{\rm case} \\ \hline
{\rm appropriate}\;{\rm value}\;{\rm for}\; H            &   H=A_0^m  &  H=A_0^{2m} \\ \hline
{\rm the}\;{\rm procedure}\;{\rm is}\;{\rm efficient}\;{\rm for}       &   n>m      &  n>2m  \\ \hline
{\rm the}\;{\rm reduced}\;{\rm bound}\;{\rm is}\;{\rm of}\;{\rm magnitude} &  A_0^\frac{m-1}{n-1} & A_0^\frac{2m-1}{n-1}\\ \hline
\end{array}
\]
This phenomenon can be detected in our examples, as well.

We start with an initial bound $A_0$ for $A$
(obtained from Theorem \ref{th1} and $\overline{|Y|}<C$)
and perform the reduction. In the following steps $A_0$
is the bound obtained in the previous reduction step. In the first reduction steps the new bound is drastically smaller than the original one.
Applying Theorem \ref{redlemma} in several steps (usually 5-10 steps) the procedure
stops by giving (almost) the same bound like the previous one. This reduced bound $A_R$
is usually between 10 and 500.

\noindent
{\bf Remark 2}\\
As it is seen we use the same lattice and our following reduction Theorem \ref{redlemma} is almost the
same as in \cite{csextic}, \cite{relthue} (see also \cite{book}).
However the approach and the way of application is completely different and
just that is our main goal. 

In the above cited papers,
following the standard arguments of \cite{relthue} we write 
$X-\alpha Y$ as a product of an element of relative norm $\mu$ and a
power product of fundamental units. Using Siegel's identity, elementary
estimates and Baker's method we arrive at a linear form in the logarithms of algebraic
numbers, which is small. The coefficients of this linear form 
in the logarithms of algebraic numbers appear in the last rows of the matrix
used instead of our matrix above to define the lattice.
The variables are the exponents of the fundamental units in the representation
of $X-\alpha Y$.

In our setting the variables are the $x_j,y_j$ and their
coefficients in (\ref{redineq}) show up in the last rows of our matrix. 
This yields a much more direct way without calculating 
fundamental units, elements of given relative norm and without applying Baker's method.

\section{Enumerating the tiny values}
\label{sectenum}

As we saw in the previous section, in our statements we assume that $\overline{|Y|}$
is not very small, which is equivalent to $A$ being large enough. The
reduction process of the previous section produces a small bound for $A$.
Hence the small values of $x_1,\ldots,x_m,y_1,\ldots,y_m$ must be tested
separately. The purpose of this section is to develop an efficient
algorithm for testing all $x_1,\ldots,x_m,y_1,\ldots,y_m$ with absolute values
less that a prescribed bound $A_R$.

For a given $m$ the number of $x_i,y_i\; (1\leq i\leq m)$ with $A\leq A_R$ is
$(2A_R+1)^{2m}$ which can still be a huge number for $m\geq 3$. 
We show how to overcome this difficulty.
(Remember that in \cite{relthue} we had to use a relatively complicated
ellipsoid method to deal with the small exponenets of the fundamental units.)

\subsection{A direct procedure}
\label{direct}

The original equation (\ref{rel}) implies
\begin{equation}
\prod_{j=1}^n(X^{(h)}-\alpha^{(hj)} Y^{(h)})=\mu^{(h)} 
\label{prod}
\end{equation}
for $h=1,\ldots, m$. We let $y_1,\ldots,y_m$ run between $-A_R$ and $A_R$ (for even $n$
it is sufficient to take nonnegative $y_1$). For each $y_1,\ldots,y_m$ we calculate 
$Y^{(h)}=y_1+\omega_2^{(h)}y_2+\ldots +\omega_m^{(h)}y_m$ and calculate the roots 
$X^{(h1)},\ldots,X^{(hn)}$ of polynomial equation (\ref{prod}) in $X^{(h)}$. 
Since we can not know which root corresponds to which conjugate, to get all solutions,
we have to solve the system of equations
\begin{eqnarray*}
x_1+\omega_2^{(1)}x_2+\ldots+\omega_m^{(1)}x_m&=&X^{(1\;i_1)}\\
&\vdots&\\
x_1+\omega_2^{(m)}x_2+\ldots+\omega_m^{(m)}x_m&=&X^{(m\; i_m)}\\\\
\end{eqnarray*}
for all possible permutations $(i_1,\ldots,i_m)$ of $(1,\ldots,m)$
and check whether the solution vector $(x_1,\ldots,x_m)$ has integer components. 

For totally real $M$ we exclude complex values of $X^{(hj)}$. 

This procedure yields $(2\cdot A_R+1)^m \; n^m$ tests. This can be used for $m=2$ but not in 
this simple form for $m\geq 3$.

\subsection{The case $m\geq 3$}
\label{chol}

For $m\geq 3$ we proceed as follows. We take a rather small initial value $A_I$ (say 10 or 20)
such that the above direct procedure of subsection \ref{direct}
can be performed to test $y_1,\ldots,y_m$ with
absolute values $\leq A_I$ within feasible CPU time.
Then we only have to consider values of $y_1,\ldots,y_m$ with $\max |y_j|>A_I$
yielding $A=\max (\max |x_j|, \max |y_j|)>A_I$.

In some steps we construct intervals $[A_s,A_S]$ the union of which covers the 
whole interval $[A_I,A_R]$. This means that in the first step we take $A_s=A_I$ and 
an $A_S$ with $A_s\leq A_S\leq A_R$. In the following step we set $A_s$ to be the former $A_S$ and take a new $A_S$, etc.

We describe now an efficient method to enumerate the variables with $A_s\leq A\leq A_S$.
For a given $h$ and $i$ specified in the reduction procedure ($1\leq h\leq m,1\leq i\leq n$) 
by (\ref{redineq}) we have
\[
|x_1+\omega_2^{(h)}x_2+\ldots + \omega_m^{(h)}x_m 
-\alpha^{(hi)}y_1-\alpha^{(hi)}\omega_2^{(h)}y_2-\ldots - \alpha^{(hi)}\omega_m^{(h)}y_m |\leq 
\]
\[
\leq d_{hi} A^{1-n} \leq d_{hi} A_s^{1-n}.
\]
We take $H=A_S\cdot A_s^{n-1}/d_{hi}$, then 
\[
H\cdot |x_1+\omega_2^{(h)}x_2+\ldots + \omega_m^{(h)}x_m 
-\alpha^{(hi)}y_1-\alpha^{(hi)}\omega_2^{(h)}y_2-\ldots - \alpha^{(hi)}\omega_m^{(h)}y_m |\leq
\]
\begin{equation}
\leq A_S
\label{redineq2}
\end{equation}

Denote by $e_1,\ldots,e_m,f_1,\ldots,f_m$ the colums of the matrix defining the lattice $\cal L$
in the preceeding section. Using the above $H$ implies that all coordinates of 
$x_1e_1+\ldots +x_me_m+y_1f_1+\dots+_mf_m$ are less than or equal to $A_S$ yielding
\begin{equation}
|x_1e_1+\ldots +x_me_m+y_1f_1+\dots+y_mf_m|^2\leq (n+1)A_S^2.
\label{norma}
\end{equation}
(In the totally real case we omit the last row.)

This defines an ellipsoid. The integer points can be enumerated by using the 
Cholesky decomposition (see M.Pohst \cite{dmv}, M.Pohst and H.Zassenhaus \cite{pz}). 
This means to construct an upper triangular matrix $R=(r_{ij})$ 
with positive diagonal entries, such that
the symmetric matrix of the above quadratic form is written as $R^T R$, that is 
(denoting here $y_1,\ldots,y_m$ by $x_{m+1},\ldots,x_{2m}$, for simplicity)
(\ref{norma}) gets the form
\[
\sum_{i=1}^{2m}\left( r_{ii}x_i+\sum_{j=i+1}^{2m}r_{ij}x_j   \right)^2\leq (n+1)A_S^2.
\]
By enumerating $x_{2m},x_{2m-1},\ldots $ etc. we also use that fact that 
\[
-A_S\leq x_i \leq A_S \;\; (1\leq i\leq 2m).
\]
Note that the Cholesky decomposition can be improved by using the Fincke-Pohst method
\cite{fipo} (see also M.Pohst \cite{dmv}, M.Pohst and H.Zassenhaus \cite{pz}),
involving LLL reduction, but then we loose the above bounds for the $x_i$.

\section{The complete algorithm}

In this section we construct the algorithm using the components of the preceeding sections.

\noindent
{\bf Problem}. Determine all solutions $X,Y\in\Z_M$ with $\overline{|Y|}<C$ of the equation
\[
N_{K/M}(X-\alpha Y)=\mu.
\]
We assume that $n>m$ if $K$ is not totally real and $n>2m$ if $K$ is totally real.

\vspace{0.5cm}

\noindent
{\bf Step 1.} Calculate the constants $c_{hji},c_{hi}\; (1\leq h\leq m,1\leq i,j\leq n,i\neq j)$ and
$c_1,c_2,c_3,c_4,c_5$ and $d_{hi}\; (1\leq h\leq m,1\leq i\leq n)$
of Section \ref{sectel}.

\vspace{0.5cm}

\noindent
{\bf Step 2.} Set $A_B=c_5\cdot \max (C,c_4) $. (This is the initial upper bound for $A$.)

\vspace{0.5cm}

\noindent
{\bf Step 3.} For $h\in\{1,\ldots,m\}$ and $i\in\{1,\ldots,n\}$ perform the reduction procedure 
of Section \ref{sectred}.\\
In the first step take $A_0=A_B$, choose a suitable $H$, perform the LLL basis reduction and calculate
the reduced bound $A_1$ by (\ref{rrr}). In the next step take $A_0=A_1$ and perform the reduction again.
Continue until the reduced bound is not any more considerably less then the previous bound.
Denote by $A_{R,h,i}$ the final reduced bound.

\vspace{0.5cm}

\noindent
{\bf Step 4.} Set $A_R=\max_{h,i} A_{R,h,i}$.

\vspace{0.5cm}

\noindent
{\bf Step 5.} If $c_5c_4> A_R$ then set $A_R=c_5c_4$. (Our arguments are only valid for 
$\overline{|Y|}>c_4$ that is $A>c_5c_4$).

\vspace{0.5cm}

\noindent
{\bf Step 6.} Enumerate the tiny values of $x_1,\ldots,x_m,y_1,\ldots,y_m$ with $A\leq A_R$
using the procedure described in Section \ref{sectenum}.
Test all possible vectors by substituting it into the equation.

\section{Computational aspects}

All our algorithms were developed in Maple \cite{maple} under Linux and the execution times 
of our examples refer to a middle category laptop. However, especially at the final enumeration,
to find appropriate values $A_s,A_S$ we made several test runs on the supercomputer
located in Debrecen, Hungary. The HPC running times were 20-50 percent shorter even
on a single node.

\section{Examples}

In this section we demonstrate our method with four explicit examples.
To have a comparison, the complete resolution of a quartic Thue equation
over a quadratic field in \cite{relthue} took about one hour.
Here we calculate "small" solutions of relative Thue
equations of degrees 6, 9, 21, over quadratic fields within a few minutes. 
Our last example is a sextic relative Thue equation over a cubic
field which could hardly be dealt with using the methods of  \cite{relthue}
(the underlying number field is of degree 18).

\subsection{Example 1}

Let $M=\Q(i)$ with integral basis $\{1,i\}$. Let $\alpha$ be a root of $f(x)=x^6+x+1$
and let $K=M(\alpha)$. Determine all $X,Y\in\Z_M$ with $\overline{|Y|}<10^{500}$
satisfying
\begin{equation}
N_{K/M}(X-\alpha Y)=X^6+XY^5+Y^6=1.
\label{p1}
\end{equation}
We get $A\leq 0.2252\cdot 10^{501}=A_B$. The reduction process ran as follows:

\[
\begin{array}{|c|c|c|c|c|c|c|}
\hline
step& A_0                 &    H      &   ||b_1||\geq        &  Digits   &    new A_0           & CPU \; time \\ \hline
1.  &0.2252\cdot 10^{501} & 10^{1003} &0.1424\cdot 10^{501}  &  1150     & 0.9637\cdot 10^{101} & 150 sec   \\ \hline 
2.  &0.9637\cdot 10^{101} & 10^{205 } &0.6095\cdot 10^{102}  &  250      & 0.1809\cdot 10^{22}  & 7 sec   \\ \hline 
3.  &0.1809\cdot 10^{22} & 10^{45} &0.1144\cdot 10^{23}  &  70     & 159562 & 2 sec   \\ \hline 
4.  & 159562 & 10^{13} &0.1009\cdot 10^{7}  &  30     & 103 & 1 sec   \\ \hline 
5.  & 103 & 10^{6} & 651.4291  &  20     & 17 & 1 sec   \\ \hline 
6.  & 17 & 10^{5} & 107.5174  &  20     & 16 & 1 sec   \\ \hline 
\end{array}
\]
The direct method of subsection (\ref{direct}) to enumerate the variables 
with absolute values $\leq A_R=16$ took 10 seconds.
Finally all solutions (up to sign) are 
\[
(x_1,x_2,y_1,y_2)=(1,0,0,0),(1,0,-1,0),(0,0,1,0).
\]

\subsection{Example 2}

Let $M=\Q(i\sqrt{2})$ with integral basis $\{1,i\sqrt{2}\}$. 
Let $\alpha$ be a root of $f(x)=x^{21}-x-1$
and let $K=M(\alpha)$. Determine all $X,Y\in\Z_M$ with $\overline{|Y|}<10^{500}$
satisfying
\begin{equation}
N_{K/M}(X-\alpha Y)=X^{21}-XY^{20}-Y^{21}=1.
\label{p2}
\end{equation}
We get $A\leq 0.2068\cdot 10^{501}=A_B$. The reduction process ran as follows:

\[
\begin{array}{|c|c|c|c|c|c|c|}
\hline
step& A_0                 &    H      &   ||b_1||\geq        &  Digits   &    new A_0           & CPU \; time \\ \hline
1.  &0.2068\cdot 10^{501} & 10^{1005} &0.1307\cdot 10^{502}  &  1200     & 0.5897\cdot 10^{26} & 420 sec   \\ \hline 
2.  &0.5897\cdot 10^{26} & 10^{55} &0.3730\cdot 10^{27}  &  80      & 99  & 5 sec   \\ \hline 
3.  & 99  & 10^{7} &  626.1309  &  20    & 6 & 2 sec   \\ \hline 
\end{array}
\]
For a larger $n$ the reduction is very efficient.
In this example $\overline{|Y|}>c_4$ implies that the main arguments 
are only valid for $\overline{|Y|}>9.9271$, that is the variables with $A\leq 20$
must be considered separately. Therefore we let the direct method 
of subsection (\ref{direct}) run for $A_R=20$. This took about 5 minutes
and resulted the solutions\\
 $(x_1,x_2,y_1,y_2)=(1,0,0,0),(1,0,-1,0),(0,0,-1,0),(-1,0,-1,0)$.

\subsection{Example 3}

This is an example for the totally real case. 
Let $M=\Q(\sqrt{2})$ with integral basis $\{1,\sqrt{2}\}$. 
Let $\alpha$ be a root of $f(x)=x^9+3x^8-5x^7+17x^6+7x^5-30x^4-x^3+16x^2-2x-1$.
This totally real nonic polynomial is taken from J.Voight \cite{vv}.
Let $K=M(\alpha)$. Determie all $X,Y\in\Z_M$ with $\overline{|Y|}<10^{500}$
satisfying
\[
N_{K/M}(X-\alpha Y)=X^9+3X^8Y-5X^7Y^2+17X^6Y^3+7X^5Y^4
\]
\begin{equation}
-30X^4Y^5-X^3Y^6+16X^2Y^7-2XY^8-Y^9=1
\label{p3}
\end{equation}
We get $A\leq 0.5379\cdot 10^{501}=A_B$. The reduction process ran as follows:

\[
\begin{array}{|c|c|c|c|c|c|c|}
\hline
step& A_0                 &    H      &   ||b_1||\geq        &  Digits   &    new A_0           & CPU \; time \\ \hline
1.  &0.5379\cdot 10^{501} & 10^{2007} &0.3402\cdot 10^{502}  &  2200     & 0.8862\cdot 10^{189} & 840 sec   \\ \hline 
2.  &0.8862\cdot 10^{189} & 10^{800} &0.5604\cdot 10^{190}  &  900     & 0.1110\cdot 10^{78} & 120 sec   \\ \hline 
3.  &0.1110\cdot 10^{78} & 10^{313} &0.7022\cdot 10^{78}  &  400     & 0.1439\cdot 10^{31} & 60 sec   \\ \hline 
4.  &0.1439\cdot 10^{31} & 10^{125} &0.9103\cdot 10^{31}  &  200     & 0.3304\cdot 10^{13} & 30 sec   \\ \hline 
5.  &0.3304\cdot 10^{13} & 10^{55} &0.2089\cdot 10^{14}  &  80     & 941870 & 30 sec   \\ \hline 
6.  & 941870 & 10^{28} &0.5956\cdot 10^{7}  &  80     & 2612 & 30 sec   \\ \hline 
7.  & 2612 & 10^{18} & 16519.5940  &  80     & 306 & 30 sec   \\ \hline 
8.  & 306 & 10^{14} & 1935.2970  &  80     & 126 & 30 sec   \\ \hline 
\end{array}
\]
The direct method 
of subsection (\ref{direct}) with $A_R=126$ executed 21 minutes. 
The solutions are\\
 $(x_1,x_2,y_1,y_2)=(1,0,0,0),(1,0,-1,0),(0,-1,-1,0),(0,0,-1,0),$\\
 $(-1,0,-1,0),(0,1,-1,0)$.

\subsection{Example 4}

Our last example demonstrates an equation with a cubic base field, $m=3$.
Let $M=\Q(\rho)$, where $\rho$ is defined by the (totally real) polynomial $x^3-x^2-3x+1$.
The field $M$ has integral basis $\{1,\rho,\rho^2\}$.
Let $\alpha$ be a root of the (totally complex) polynomial $f(x)=x^6+2x^5+3x^4+21$.
Let $K=M(\alpha)$. Determine all $X,Y\in\Z_M$ with $\overline{|Y|}<10^{500}$
satisfying
\begin{equation}
N_{K/M}(X-\alpha Y)=X^6+2X^5Y+3X^4Y^2+21Y^6=1.
\label{p4}
\end{equation}
We get $A\leq 0.4268\cdot 10^{501}=A_B$. 
In this example we used the constants $c_{hi}$ and $d_{hi}$
calculated for the given case $h,i$ (while in all other examples
we used values valid for all cases). This resulted a better reduction,
giving a reduced bound about 15 procent sharper.

The reduction process ran as follows:

\[
\begin{array}{|c|c|c|c|c|c|c|}
\hline
step& A_0                 &    H      &   ||b_1||\geq        &  Digits   &    new A_0           & CPU \; time \\ \hline
1.  &0.4268\cdot 10^{501} & 10^{1513} &0.1277\cdot 10^{503}  &  1700     & 0.9649\cdot 10^{203} & 600 sec   \\ \hline 
2.  &0.9649\cdot 10^{203} & 10^{615} &0.2888\cdot 10^{205}  &  800     & 0.8196\cdot 10^{83} & 160 sec   \\ \hline 
3.  &0.8196\cdot 10^{83} & 10^{255} &0.2453\cdot 10^{85}  &  350     & 0.8465\cdot 10^{35} & 60 sec   \\ \hline 
4.  &0.8465\cdot 10^{35} & 10^{111} &0.2534\cdot 10^{37}  &  200     & 0.5308\cdot 10^{16} & 30 sec   \\ \hline 
5.  &0.5308\cdot 10^{16} & 10^{53} &0.1589\cdot 10^{18}  &  100     &  0.9237\cdot 10^{8} & 20 sec   \\ \hline 
6.  &0.9237\cdot 10^{8} & 10^{29} &0.2764\cdot 10^{10}  &  60     &  52170 & 20 sec   \\ \hline 
7.  & 52170 & 10^{19} &0.1561\cdot 10^{7}  &  50     &  2328 & 20 sec   \\ \hline 
8.  & 2328 & 5\cdot 10^{14} &  69684.4896  &  50     &  598 & 20 sec   \\ \hline 
9.  & 598 & 8\cdot 10^{12} &  17900.0536  &  50     &  343 & 20 sec   \\ \hline 
10.  & 343 & 4\cdot 10^{12} &  10267.0876  &  50     &  334 & 20 sec   \\ \hline 
\end{array}
\]

Note that in our example the direct method of Subsection \ref{direct} with $A_R=334$ yields to
test 
\[
(2\cdot A_R+1)^m n^m=(2\cdot 334 + 1)^3\cdot 6^3=64.674.354.744
\]
cases which is completely impossible.

We executed the direct method of Subsection \ref{direct} with $A_I=10$
taking 7 minutes. To cover the interval [10,334] we applied the algorithm of Subsection
\ref{chol} in several steps using 200 digits accuracy.
It is worthy to choose $A_s$ and $A_S$ with a relatively large difference
so that $H$ also becomes large and inequality (\ref{redineq2}) becomes
an efficient filter.
We performed four steps of the algorithm of subsection \ref{chol}:
\[
\begin{array}{|c|c|c|c|c|c|c|}
\hline
step& A_s     &  A_S     & CPU \; time \\ \hline
1.  & 10  &  50  &  3 \; min \\ \hline
2.  & 50  &  100  &  4 \; min \\ \hline
3.  & 100  &  150  &  10 \; min \\ \hline
4.  & 100  &  334  &  40 \; min \\ \hline
\end{array}
\]
We had several test runs showing an optimal segmentation of the interval [10,334].
It turned out that the running time for the last step with $A_S=334$ is not significantly
different with $A_s=100,150,200,250$. Therefore it was not worthy to split this
interval into further parts.

The only solution (up to sign) of equation (\ref{p4}) is $(x_1,x_2,y_1,y_2)=(1,0,0,0)$.

\end{document}